\documentclass[reqno,12pt]{article}
\usepackage{amsfonts}  \usepackage{amssymb}\usepackage{pstricks} 

\def\alert#1{\smallskip{\hskip\parindent\vrule%
\vbox{\advance\hsize-2\parindent\hrule\smallskip\parindent.4\parindent%
\narrower\noindent#1\smallskip\hrule}\vrule\hfill}\smallskip}

\newcommand{\CCC}{\mathcal{C}} \newcommand{\DDD}{\mathcal{D}}
\newcommand{\FFF}{\mathcal{F}}       \newcommand{\GGG}{\mathcal{G}}  \newcommand{\III}{\mathcal{I}}

\newcommand{\MMM}{\mathcal{M}} 
\newcommand{\OOO}{\mathcal{O}}
    \newcommand{\QQQ}{\mathcal{Q}}  \newcommand{\RRR}{\mathcal{R}}
\newcommand{\VVV}{\mathcal{V}}

\newcommand{\riff}{\Longrightarrow}

       \newcommand{\NN}{\mathbb{N}}
\newcommand{\xx}{\mathbf{x}}  
\newcommand{\vv}{\mathbf{v}}\newcommand{\yy}{\mathbf{y}}

\newtheorem{thm}{Theorem}[section]
\newtheorem{lem}[thm]{Lemma}
\newtheorem{cor}[thm]{Corollary}

\newtheorem{que}[thm]{Question}
\newtheorem{define}[thm]{Definition}

\usepackage{amssymb}
\usepackage{amsfonts}
\def\alert#1{\smallskip{\hskip\parindent\vrule%
\vbox{\advance\hsize-2\parindent\hrule\smallskip\parindent.4\parindent%
\narrower\noindent#1\smallskip\hrule}\vrule\hfill}\smallskip}

\newcommand{\KK}{\mathbb{K}}

\newcommand{\JJ}{\mathbb{J}}

\newcommand{\RR}{\mathbb{R}}

\newcommand{\ifff}{\Longleftrightarrow}   

\usepackage{amssymb}
\usepackage{amsfonts}
\def\alert#1{\smallskip{\hskip\parindent\vrule%
\vbox{\advance\hsize-2\parindent\hrule\smallskip\parindent.4\parindent%
\narrower\noindent#1\smallskip\hrule}\vrule\hfill}\smallskip}

\usepackage{pstricks}

\begin{document}

{\large \bf 
\begin{flushright} {\bf \large Pre-kites:}\\
{\bf \large Simplices having a regular facet}
\end{flushright}}

\vspace{.4cm}
\begin{flushright}
Mowaffaq Hajja\\
Department of Mathematics\\
Yarmouk University\\
Irbid -- Jordan\\
\vspace{.3cm}
\noindent mowhajja@yahoo.com\\~\\

\vspace{.4cm}
Mostafa Hayajneh\\
Department of Mathematics\\
Yarmouk University\\
Irbid -- Jordan\\
\vspace{.3cm}
\noindent hayaj86@yahoo.com\\~\\

\vspace{.3cm}

Ismail Hammoudeh\\
Department of Physics\\
Amman National  University\\
Salt -- Jordan\\
\vspace{.3cm}
\noindent ihammoudeh@ammanu.edu.jo\\
\noindent ismaa3iil@gmail.com\\~\\
\end{flushright}

\vspace{.4cm}
\begin{quote}
\small {\bf Abstract.}  The  investigation of the relation  among the distances of an arbitrary point in the Euclidean space $\RR^n$ to the vertices of a regular $n$-simplex in that space has led us  to the study of simplices having a regular facet. Calling an $n$-simplex with a regular facet an $n$-pre-kite,
we  investigate, in the spirit of \cite{ET}, \cite{BAG-1}, \cite{RM-ortho}, and \cite{RM-balloon}, and using tools from linear algebra, the degree of regularity implied by the coincidence of any two of the classical centers of such simplices.
We also prove that if $n \ge 3$,  then the intersection of the family of $n$-pre-kites  with any of the four known special families
is the family of $n$-kites, thus extending  the result  in \cite{kites}.
A basic tool is a closed form of a determinant that arises in the context of a certain Cayley-Menger determinant, and that generalizes several determinants that appear in \cite{RM-ortho},  \cite{RM-balloon}, and \cite{impurity}.  
Thus the paper is a further testimony to the special role that linear algebra plays in higher dimensional geometry.
\end{quote}

\vspace{-.25cm}
\begin{quote}
{\small {\bf Mathematics Subject Classification (2010).} Primary 52B11; Secondary 52B12, 52B15, 51M20, 52B10.}
\end{quote}

\vspace{-.25cm}
\begin{quote}
{\small {\bf Keywords.}  affine hull, affine independence, Cayley-Menger determinant, centroid, cevian, circumcenter,     circumscriptible simplex, incenter, inner Cayley-Menger determinant,  isodynamic simplex, isogonic simplex, equiareal simplex, equiradial simplex, kite, orthocentric simplex, orthogonal complement,  Pompeiu's theorem, pre-kite, special tetrahedra, tetra-isogonic simplex, well distributed edge lengths}
\end{quote}

\section{Introduction} \label{Int}
The distances $t_1,\cdots,t_{n+1}$ between the vertices of a regular $n$-simplex $S$ of edge length $t_0$ and  an arbitrary point $P$ in its affine hull are related by the elegant relation    
\begin{eqnarray} \label{REL}
(n+1) \sum_{j=0}^{n+1} t_j^4 &=& \left( \sum_{j=0}^{n+1} t_j^2 \right)^2; 
\end{eqnarray}
see \cite{Bentin-1} for a very short proof, and see \cite{HHSh} for a proof that (\ref{REL}) is essentially the only relation that exists among the quantities $t_0,\cdots, t_{n+1}$, for  fixed $t_0$. Figure 1 below illustrates the case when $d=2$, i.e., when the regular simplex is an equilateral triangle.  
\begin{center}
~~\psset{unit=.3cm}
\begin{pspicture}(-11,-5)(11,20)
\psline[linestyle=dashed](-15,8)(-10,0)
\psline[linestyle=dashed](-15,8)(10,0)
\psline[linestyle=dashed](-15,8)(0,17.32)
\rput(-13.5,4){\small $t_2$} \rput(-2.5,3){\small $t_3$} \rput(-8.5,13.1){\small $t_1$}



\psline(10,0)(0,17.32)
\psline(-10,0)(0,17.32)
\psline(-10,0)(10,0)

\rput(0,-1){\small $t_0$} \rput(-6,8.7){\small $t_0$} \rput(6,8.7){\small $t_0$}
\rput(0,18){\small $\vv_1$} \rput(-10.5,-1){\small $\vv_2$} \rput(10.5,-1){\small $\vv_3$}
\rput(-16,8){\small $P$}

\rput(0,-2.5){\small \bf Figure 1}
\rput(0,-4.5){\small \bf Illustrating (\ref{REL}) in the case $d=2$}
\end{pspicture}
\end{center}
 
The relation  (\ref{REL}) has been a source of fascination, and its special case $n=2$ has been a source of  inspiration to problem lovers. Problems that refer to  Figure 1 and that give numerical values of three of the variables $t_0,t_1,t_2,t_3$ and ask about the fourth variable  can, in the absence of the relation (\ref{REL}), be thought of as challenging. Such problems have appeared in \cite{500}, \cite{Graham}, \cite{Wagon}, \cite{CMJ}, and possibly others. One of the features of Figure 1 is the fact that the lengths $t_1, t_2, t_3$ can serve as the side lengths of a triangle, i.e., they satisfy the triangle inequality. This non-obvious and interesting fact is attributed to the Romanian mathematician Dimitrie Pompeiu (1873--1954), and now carries his name. It, too, appeared, with different proofs quite frequently; see, for example,  \cite{Titu}, \cite{Titu-2}, \cite{Fomin}, \cite{Alsina}, \cite{mini}, \cite{Delights}, \cite{Roy}, and \cite{Putnam}.

Moving the point $P$ outside the affine hull of the regular $n$-simplex $S$ results in an $(n+1)$-simplex $(S,P)$ having  $S$ as a facet. This led us to consider such simplices. Thus we call an $(n+1)$-simplex having a regular facet  an $(n+1)$-{\it pre-kite}, and we investigate their properties. When the edges emanating from the point $P$ to the vertices of $S$ are equal, then the pre-kite $(S,P)$ is what was called a {\it kite} in  \cite{RM-ortho} and in \cite{kites}.

In this paper, we find a closed form of the Cayley-Menger determinant of a pre-kite. A main determinant that comes up generalizes several determinants that have appeared in \cite{RM-ortho}, \cite{RM-balloon}, and \cite{impurity},
and possibly other places. By setting the Cayley-Menger  determinant $M$ equal to 0, we obtain  the relation (\ref{REL}) mentioned above, and we use this relation to give a new proof of Pompeiu's theorem. We also investigate the degree of regularity implied by the coincidence of two of the classical centers of a pre-kite, and we find the  intersection of the family of pre-kites with any of the known special families of orthocentric, circumscriptible, isodynamic, and tetra-isogonic simplices.

\bigskip
The paper is organized as follows. In Section 2, we define 
the Cayley-Menger and the inner Cayley-Menger determinants $\CCC$ and $\DDD$ of a simplex $S$, and we recall the formulas that give the volume and circumradius of $S$ in terms of $\CCC$ and $\DDD$. In Section 3, we find a closed form of a certain determinant that will be used to evaluate $\CCC$ and $\DDD$. This determinant has, as special cases, several determinants that have appeared in the earlier literature. Section 4 applies the formulas in Section 3 to pre-kites to give a new derivation of (\ref{REL}) and to give a new proof of Pompeiu's theorem. In Section 5, we prove that a non-regular pre-kite cannot have more than two different  regular facets, and we also  characterize the positive numbers that can serve as edge lengths of a two-apexed pre-kite.
In Section 6, we prove that if the circumcenter and centroid of a pre-kite coincide, then it is regular. The same holds if the circumcenter and incenter coincide. We also prove that if $n \le 5$, and if the centroid and the incenter of an $n$-pre-kite coincide, then it is regular, and we provide examples of non-regular $n$-pre-kites, $n \ge 6$,  in which the centroid and incenter coincide. We also prove that if the Fermat-Torricelli point of a pre-kite coincides with either the centroid  or the circumcenter, then it is regular. In Section 7, we prove that if a pre-kite of dimension $n \ge 3$ belongs to any of the four known special families of orthocentric, circumscriptible, isodynamic, and tetra-isogonic $n$-simplices, then it is a kite.

\bigskip It is not unusual to use tools from linear algebra, such as properties
of certain matrices and determinants, in investigations pertaining to higher-dimensional geometry. The importance of such tools is manifested, for example, in the proofs of the higher dimensional analogues of the theorems of Pythagoras, as in \cite{Lin} and \cite{Bhatia}, Napoleon, as  in \cite{Weiss} and \cite{Spirova}, the law of sines, as in \cite{Leng}, coincidences of centers, as in \cite{ET}, \cite{BAG-1}, \cite{RM-ortho}, and \cite{RM-balloon}, and the open mouth theorem, as in \cite{OMT}.

\section{The Cayley-Menger determinant, and formulas for the volume and circumradius of an $n$-simplex} This section puts together known formulas for 
the volume and circumradius of an $n$-simplex in terms of the Cayley-Menger and inner Cayley-Menger determinants of $S$.

We recall that  an {\it $n$-simplex}, $n \ge 2$, is defined to be  the convex hull $S = [A_0,\cdots, A_n]$ of $n+1$ affinely independent points $A_0, \cdots, A_n$ in a Euclidean space $\RR^m$, $m \ge n$. For $0 \le j\le n$, $A_j$ is called the $j$-th {\it vertex} of $S$, and the $(n-1)$-simplex obtained from $S$ by  removing the vertex $A_j$ is called the  $j$-th {\it facet} of $S$ and is denoted by $S_j$.  The simplex obtained from $S$ by  removing any number of vertices  is called a {\it face} of $S$.

Let $S = [A_0,\cdots,A_n]$ be an $n$-simplex, and let 
\begin{eqnarray} \label{aij} \| A_i - A_j \|^2 = a_{i,j}, ~~0 \le i, j \le n.\end{eqnarray}
The  {\it Cayley-Menger determinant} $\CCC = \CCC (S)$ of $S$ is the $(n+2) \times (n+2)$ determinant whose entries $c_{i,j}$, $-1 \le i, j \le n$, are defined by  
\begin{eqnarray} \label{cij} c_{i,j} &=& \left\{
\begin{array}{ll}
0 & \mbox{if $i=j$},\\
1 & \mbox{if $i=-1$ and $j \ne-1$},\\
1 & \mbox{if $j=-1$ and $i \ne-1$},\\
a_{i,j} & \mbox{otherwise};
\end{array}
\right.
\end{eqnarray}
see, for example,  \cite[\S 9.7.3.1, pp.~237--238]{Berger}. Thus
\begin{eqnarray}
\CCC=\CCC (S) &=&  \left|
\begin{array}{cccccccc}
0&1&1&1&\cdots&\cdots&1&1\\
1&0&a_{0,1}&a_{0,2}&\cdots&\cdots&a_{0,n-1}&a_{0,n}\\
1&a_{1,0}&0&a_{1,2}&\cdots&\cdots&a_{1,n-1}&a_{1,n}\\
\cdots  &\cdots &\cdots  &\cdots &\cdots&\cdots &\cdots &\cdots\\
\cdots  &\cdots &\cdots  &\cdots &\cdots&\cdots &\cdots &\cdots\\
1&a_{n-1,0}&a_{n-1,1}&a_{n-1,2}&\cdots&\cdots&0&a_{n-1,n}\\
1&a_{n,0}&a_{n,1}&a_{n,2}&\cdots&\cdots&a_{n,n-1}&0\\
\end{array}
\right|. \label{CCC}
\end{eqnarray}
If $\VVV = \VVV (S)$ is the volume (i.e., the $n$-dimensional Lebesgue measure or content) of $S$, then it is well known that
\begin{eqnarray}\label{V}
(-1)^{n+1} 2^n (n!)^2  \VVV^2 &=& \CCC;
\end{eqnarray}
see, for example, \cite[(5.1), \S 5, Chapter VIII, p.~125]{Sommerville} and \cite{Ivanoff}. 

The determinant obtained from $\CCC$ by deleting the uppermost row and the leftmost column will be  denoted  by $\DDD = \DDD (S)$, and will be referred to as the {\it inner Cayley-Menger determinant} of $S$.
Thus 
\begin{eqnarray}
\DDD = \DDD (S) &=&  \left|
\begin{array}{ccccccc}
0&a_{0,1}&a_{0,2}&\cdots&\cdots&a_{0,n-1}&a_{0,n}\\
a_{1,0}&0&a_{1,2}&\cdots&\cdots&a_{1,n-1}&a_{1,n}\\
\cdots &\cdots  &\cdots &\cdots&\cdots &\cdots &\cdots\\
\cdots &\cdots  &\cdots &\cdots&\cdots &\cdots &\cdots\\
a_{n-1,0}&a_{n-1,1}&a_{n-1,2}&\cdots&\cdots&0&a_{n-1,n}\\
a_{n,0}&a_{n,1}&a_{n,2}&\cdots&\cdots&a_{n,n-1}&0\\
\end{array}
\right|. \label{DDD}
\end{eqnarray}

The determinants $\CCC = \CCC (S)$ and $\DDD = \DDD (S)$ are used in \cite{Ivanoff} to express  the circumradius $\RRR = \RRR (S)$ of $S$ as
\begin{eqnarray}\label{R}
\RRR^2 &=& \frac{- \DDD}{2 \CCC}. 
\end{eqnarray}

\section{A special determinant} \label{special-det}
In this section, we consider the determinant $\KK(n;z;x_1,\cdots,x_n;y_1,\cdots,y_n;a;b)$ defined by (\ref{K}) below, and we evaluate it in closed form in Theorem \ref{KK}. This will then be used in Theorems \ref{CjDj} and \ref{CCCDDD} to find formulas for the volumes and circumradii of pre-kites and their facets. These formulas will in turn be used to determine the degree of regularity implied by the coincidence of any two  of the classical centers of a pre-kite. Note that the special cases of $\KK(n;z;x_1,\cdots,x_n;y_1,\cdots,y_n;a;b)$ when
$[z=1,~x_i=y_i]$, $[z=0,~x_i=y_i,~a=1,~b=-1]$, and  $[z=0,~x_i=y_i]$
have appeared in \cite{RM-ortho}, \cite{RM-balloon},  and \cite{impurity}, respectively,  where they were instrumental in establishing the results there. 

 \bigskip
We start with defining the special determinants $\JJ$ and $\KK$. In all that follows, $a, b, z, x_j, y_j$ stand for real numbers for all non-negative integers $j$, but the treatment may still hold over other rings.

\begin{define} \label{DefJK}
{\em  The determinant  $\JJ (n;z;a;b)$, $n \in \NN$, is  the $n \times n$  determinant that has $b$ on every entry on the main diagonal and $a$ everywhere else. The determinant  $\KK (n;z;x_1,\cdots,x_n;y_1,\cdots,y_n;a;b)$ is  the $(n+1) \times (n+1)$ determinant $\left(d_{i,j}\right)_{0 \le i, j \le n}$ whose 0-th row is $[z,y_1,\cdots,y_n]$,  whose 0-th column is $[z,x_1,\cdots,x_n]^t$, and whose subdeterminant $\left(d_{i,j}\right)_{1 \le i, j \le n}$ is $\JJ (n;a;b)$. More formally, the entries $d_{i,j}, 0 \le i, j \le n,$ are given by
\begin{eqnarray}
d_{i,j} &=& \left\{
\begin{array}{cl}
z & \mbox{if $i=j=0$},\\
x_i & \mbox{if $j=0$ and $1 \le i \le n$},\\
y_j & \mbox{if $i=0$ and $1 \le j \le n$},\\
b & \mbox{if  $1 \le i = j \le n$},\\
a & \mbox{otherwise}.
\end{array} \right. \label{K-formal}
\end{eqnarray}
Thus
\begin{eqnarray}
\JJ (n;a;b) &=& \left|
\begin{array}{ccccc}
b&a&a&\cdots&a\\a&b&a&\cdots&a\\a&a&b&\cdots&a\\
\cdots &\cdots &\cdots&\cdots &\cdots \\
a&a&a&\cdots&b
\end{array}
\right|,~~~\mbox{(of size $n\times n$)} \nonumber\\
&=& ((n-1)a+b) (b-a)^{n-1}, \mbox{~by Lemma 3.1 of \cite{impurity}}. \label{J}
\end{eqnarray}
The formula above has also appeared as Lemma 7.1 in \cite{RM-balloon} and as Lemma 3.6 in \cite{RM-ortho}. Also,
\begin{eqnarray}
\KK (n;z;\xx;\yy;a;b) &=&\KK (n;z;x_1,\cdots,x_n;y_1,\cdots,y_n;a;b) \nonumber\\ && \nonumber\\
&=& \left|
\begin{array}{ccccccc}
z&y_1&y_2&\cdots&y_j&\cdots&y_n\\
x_1&b&a&\cdots&a&\cdots&a\\
x_2&a&b&\cdots&a&\cdots&a\\
\cdots  &\cdots &\cdots  &\cdots &\cdots&\cdots &\cdots \\
x_j&a&a&\cdots&b&\cdots&a\\
\cdots  &\cdots &\cdots  &\cdots &\cdots&\cdots &\cdots \\
x_n&a&a&\cdots&a&\cdots&b
\end{array}
\right|. \label{K} \end{eqnarray}
}\end{define}

\begin{thm} \label{KK}
Let $n \ge 1$, and let
\begin{eqnarray} \label{xxyy}
\xx  = (x_1,\cdots,x_n),~~\yy = (y_1,\cdots, y_n),
 \end{eqnarray}
and 
\begin{eqnarray} \label{MxMyMxy}
\MMM_{\xx} = \sum_{j=1}^n  x_j,~~\MMM_{\yy} = \sum_{j=1}^n  y_j,~~
\MMM_{\xx\yy} = \sum_{j=1}^n  x_jy_j. \end{eqnarray}
Let  $\KK = \KK (n;z;\xx;\yy;a;b)$ be the determinant defined by (\ref{K}).  
Then
\begin{eqnarray}
&& \KK (n;z;\xx;\yy;a;b) \nonumber\\
&=& (b-a)^{n-2} \left[
((n-1)a+b)
\left(z(b-a) - \MMM_{\xx\yy}\right) + a \MMM_{\xx} \MMM_{\yy}\right].
\label{K=} \end{eqnarray}
\end{thm}

\vspace{.0cm} \noindent {\it Proof.} 
We proceed by induction. For $n=1$, the statement is trivial, being nothing but
\begin{eqnarray*}\KK (1;z;x_1;y_1;a;b) &=& \left| \begin{array}{cc} z&y_1\\x_1&b \end{array} \right|= zb-x_1y_1.
\end{eqnarray*}
Suppose now that (\ref{K=}) holds for $n= r$ for some $r \ge 1$. We are to show that it holds for $n = r+1$. Thus we let
\begin{eqnarray} 
\xx' &=& \left(x_1,\cdots,x_{r+1}\right),~~
\yy' ~=~ \left(y_1,\cdots,y_{r+1}\right),\label{x'}\\
\MMM_{\xx}' &=& \sum_{j=1}^{r+1} x_j,~~\MMM_{\yy}' ~=~ \sum_{j=1}^{r+1} y_j,~~
\MMM_{\xx\yy}' ~=~ \sum_{j=1}^{r+1}  x_jy_j,\label{M'} \end{eqnarray}
and we show that
\begin{eqnarray}
&&\KK(r+1;z;\xx';\yy';a;b) \nonumber \\
&=& (b-a)^{r-1}\left[
  (ra+b) \left(z(b-a) - \MMM_{\xx\yy}'\right) + a \MMM_{\xx}'
  \MMM_{\yy}' \right]. \label{M-r+1}
\end{eqnarray}
For simplicity, we denote $\KK(r+1;z;\xx',\yy';a;b)$ by $K$, and we refer to its rows and columns as the 0-th, the first, etc. Thus the   0-th row of $K$ is $[z,y_1,\cdots,y_{r+1}]$, and the 0-th column is $[z,x_1,\cdots,x_{r+1}]^t$. Expanding $K$ along the 0-th row, we obtain
\begin{eqnarray}
K &=& z C_0 +\sum_{j=1}^{r+1} (-1)^{j} y_j C_j, \label{C0Cj}
\end{eqnarray}
where $C_j$ is the $(0,j)$-th minor of $K$. 
Since $C_0$ is the case $n = r+1$ of the determinant given in (\ref{J}), it is clear that
\begin{eqnarray} \label{C0-1}
 C_0 &=& \JJ (r+1;a;b) ~=~ 
(ra+b)(b-a)^r, \mbox{~by (\ref{J})}.
\end{eqnarray} 
 To  calculate $C_j$, $1 \le j \le r+1$, we recall that $C_j$ is obtained from $K$ by deleting the 0-th row and $j$-th column, and we let $R_1, \cdots, R_{r+1}$ be the rows of $C_j$. Notice that $R_1=[x_1,b,a,\cdots,a]$ and $R_j = [x_j,a,a,\cdots,a]$, for $j > 1$. Let $E_j$ be the determinant obtained from $C_j$ by moving  $R_j$ to the very top (with $E_1=C_1$). Thus the rows of $E_j$ are $R_j,R_1,\cdots,R_{j-1},R_{j+1},\cdots,R_{r+1}$, i.e.,  $R_{\sigma (1)},\cdots,R_{\sigma(r+1)}$, where $\sigma$ is the cyclic permutation $(1~~j~~j-1~~j-2~~\cdots~~2)$. Thus
$C_j = (-1)^{j-1} E_j$. Also, the uppermost row of $E_j$ is $[x_j,a,\cdots,a]$, the leftmost column is $[x_j, x_1,\cdots,x_{j-1},x_{j+1},\cdots,x_{r+1}]^t$, and the remaining entries form $\JJ(r;a;b)$. 
Therefore 
\begin{eqnarray}
C_j &=& (-1)^{j-1} E_j \nonumber\\
&=& (-1)^{j-1} \KK(r;x_j;x_1,\cdots,x_{j-1},x_{j+1},\cdots,x_{r+1};a,\cdots,a;a;b)\nonumber\\
&=& (-1)^{j-1} (b-a)^{r-2} \left[((r-1)a+b)(x_j (b-a)-a(\MMM_{\xx}' - x_j)) + a (\MMM_{\xx}' - x_j)(ra)\right]\nonumber \\ 
&=& (-1)^{j-1} (b-a)^{r-2} \left[((r-1)a+b)x_j (b-a)-a(\MMM_{\xx}' - x_j) ((r-1)a+b-ra)\right] \nonumber\\
&=& (-1)^{j-1} (b-a)^{r-2+1} \left[((r-1)a+b)x_j-a(\MMM_{\xx}' - x_j)\right]
\nonumber\\
&=& (-1)^{j-1} (b-a)^{r-1} \left[(ra+b)x_j-a\MMM_{\xx}'\right].\label{Cj-1}
\end{eqnarray}
Using (\ref{C0Cj}), (\ref{C0-1}), and (\ref{Cj-1}), we obtain
\begin{eqnarray*}
K&=&z(ra+b)(b-a)^r+\sum_{j=1}^{r+1}(-1)^{j}y_j(-1)^{j-1}(b-a)^{r-1}\left[ (ra+b)x_j-a\MMM_{\xx}'\right] \\
&=&(b-a)^{r-1}\left[ z(ra+b)(b-a)-\sum_{j=1}^{r+1}y_j\left( (ra+b)x_j-a \MMM_{\xx}'\right) \right] \\
&=&(b-a)^{r-1}\left[ (ra+b)\left( z(b-a)-\sum_{j=1}^{r+1}x_jy_j\right) + a\MMM_{\xx}' \sum_{j=1}^{r+1}y_j \right] \\
&=&(b-a)^{r-1}\left[ (ra+b)\left( z(b-a)-\MMM_{\xx\yy}'\right) + a\MMM_{\xx}'\MMM_{\yy}' \right], 
\end{eqnarray*}
as desired. This completes the proof. \hfill $\Box$

\section{Pre-kites and formulas for their volumes and circumradii} \label{prekites}
In this section, we introduce the new family of pre-kites, and we use Theorem \ref{KK} to derive formulas for the volumes and circumradii of these simplices and of their facets. These formulas will be used in Section \ref{coincidence} to investigate the degree of regularity implied by the coincidence of two of the classical centers of an $n$-pre-kite $S$, $n \ge 2$.

We shall call  the $n$-simplex $S=[A_0,\cdots,A_n]$ an {\it $n$-pre-kite}  if one of the  facets $S_j$ is a regular $(n-1)$-simplex. In this case, we call $A_j$ an {\it apex}, and $S_j$ a {\it base} of $S$. Actually there will be no harm in referring to these as ``the" apex and ``the" base, although  an $n$-pre-kite can have more than one apex (and hence more than one base), as we shall see later in Section \ref{Jaws}.

Notice that if the $n$-simplex $S=[A_0,\cdots,A_n]$ is an $n$-pre-kite with apex $A_0$, and if the  lengths of the edges that emanate from $A_0$ are all equal, then $S$ is what was called an $n$-kite in \cite{kites} and in other papers. Notice also that all the facets (and hence all the faces) of a pre-kite are also pre-kites.

If  the $n$-simplex $S=[A_0,\cdots,A_n]$ is an $n$-pre-kite with apex $A_0$, and if
\begin{eqnarray} \label{v}
\|A_i - A_j\|^2  = \left\{
\begin{array}{cl}
v_i & \mbox{~~if $j=0$ and  $1 \le i \le n$},
\\
u & \mbox{~~$i \ne j$ and  $1 \le i, ~ j \le n$},
\end{array}
\right.
\end{eqnarray}
then we shall  denote $S$ by $PK[n;u;v_1,\cdots,v_n]$. 
Note that this $n$-pre-kite 
is an $n$-kite precisely when $v_1 = \cdots = v_n$.


\begin{thm} \label{CCCDDD} Let   
\begin{eqnarray} \label{vv} \vv &=& (v_1,\cdots,v_n),  \end{eqnarray}  and let $S=[A_0,\cdots,A_n]$  be the $n$-pre-kite with apex $A_0$ defined by
\begin{eqnarray} \label{pre-k-S} S&=& PK[n;u;\vv] ~=~ PK[n;u;v_1,\cdots,v_n], 
  \end{eqnarray}
where $n \ge 3$ and $u > 0$. Then the Cayley-Menger determinant of $S$ is given by
\begin{eqnarray}
\CCC (PK[n;u;\vv]) &=&  
(-u)^{n-2}  \left[n (u^2 + v_1^2 + \cdots + v_n^2) -  (u + v_1 + \cdots  + v_n)^2  \right], \label{CCC-PK-2}
\end{eqnarray}
and the inner Cayley-Menger determinant of $S$ is given by
\begin{eqnarray}
\DDD (PK[n;u;\vv]) &=&  
(-u)^{n-1} \left[
(n-1)  (v_1^2+\cdots + v_n^2) -  (v_1+\cdots+v_n)^2 \right]. 
\label{DDD-PK-2}
\end{eqnarray}
\end{thm}
\vspace{.2cm} \noindent {\it Proof.}  The Cayley-Menger determinant of the pre-kite $PK[n;u;\vv]$ is given by
\begin{eqnarray}
\CCC (PK[n;u;\vv]) &=&  \left|
\begin{array}{cccccccc}
0&1&1&1&\cdots&\cdots&1&1\\
1&0&v_1&v_{2}&\cdots&\cdots&v_{n-1}&v_{n}\\
1&v_{1}&0&u&\cdots&\cdots&u&u\\
\cdots  &\cdots &\cdots  &\cdots &\cdots&\cdots &\cdots &\cdots\\
\cdots  &\cdots &\cdots  &\cdots &\cdots&\cdots &\cdots &\cdots\\
1&v_{n-1}&u&u&\cdots&\cdots&0&u\\
1&v_{n}&u&u&\cdots&\cdots&u&0\\
\end{array}
\right|. \label{CCC-PK}
\end{eqnarray}
Multiplying the uppermost row by $u$ and interchanging it  with the next row, and then multiplying the leftmost column  by $u$ and  interchanging it with the next column, we obtain 
 \begin{eqnarray}
u^2 \CCC (PK[n;u;\vv]) &=&  \left|
\begin{array}{cccccccc}
0&u&v_1&v_{2}&\cdots&\cdots&v_{n-1}&v_{n}\\
u&0&u&u&\cdots&\cdots&u&u\\
v_1&u&0&u&\cdots&\cdots&u&u\\
\cdots  &\cdots &\cdots  &\cdots &\cdots&\cdots &\cdots &\cdots\\
\cdots  &\cdots &\cdots  &\cdots &\cdots&\cdots &\cdots &\cdots\\
v_{n-1}&u&u&u&\cdots&\cdots&0&u\\
v_{n}  &u&u&u&\cdots&\cdots&u&0\\
\end{array}
\right| \nonumber
 \\ \vspace{.08cm} \nonumber\\ 
 &=& \KK (n+1;0;u,v_1,\cdots,v_n;u,v_1,\cdots,v_n;u,0)
 \nonumber \\ \vspace{.03cm} \nonumber\\
 &=& (-u)^{n-1} \left[(-nu)(u^2 + v_1^2 + \cdots + v_n^2) + u
 (u + v_1 + \cdots  + v_n)^2\right] \nonumber\\ \vspace{.03cm} \nonumber\\
 &=&  (-u)^{n}  \left[n (u^2 + v_1^2 + \cdots + v_n^2) -  (u + v_1 + \cdots  + v_n)^2  \right]. \nonumber
\end{eqnarray}
Therefore
\begin{eqnarray*}
\CCC (PK[n;u;\vv]) &=&  
(-u)^{n-2}  \left[n (u^2 + v_1^2 + \cdots + v_n^2) -  (u + v_1 + \cdots  + v_n)^2  \right], 
\end{eqnarray*}
as desired.

We now calculate the  inner Cayley-Menger determinant $\DDD (PK[n;u;\vv])$  of the $n$-pre-kite $PK[n;u;\vv]$.
\begin{eqnarray*}
\DDD (PK[n;u;\vv]) &=&  \left|
\begin{array}{ccccccc}
0&v_{1}&v_{2}&\cdots&\cdots&v_{n-1}&v_{n}\\
v_{1}&0&u&\cdots&\cdots&u&u\\
\cdots &\cdots  &\cdots &\cdots&\cdots &\cdots &\cdots\\
\cdots &\cdots  &\cdots &\cdots&\cdots &\cdots &\cdots\\
v_{n-1}&u&u&\cdots&\cdots&0&u\\
v_{n}&u&u&\cdots&\cdots&u&0\\
\end{array}
\right| \nonumber \\ \vspace{.12cm} \nonumber \\ 
&=& \KK (n;0;v_1,\cdots,v_n;v_1,\cdots,v_n;u,0) \nonumber \\ \vspace{.05cm} 
&=& 
(-u)^{n-2} \left[
(n-1) (-u) (v_1^2+\cdots + v_n^2) +u  (v_1+\cdots+v_n)^2 \right] \nonumber
\\ &=&
(-u)^{n-1} \left[
(n-1)  (v_1^2+\cdots + v_n^2) -  (v_1+\cdots+v_n)^2 \right]. 
\end{eqnarray*}

This completes the proof. \hfill $\Box$

\bigskip
The next theorem is immediate, but we record it for ease of reference.

\begin{thm} \label{CjDj} Let $n \ge 3$, and let  $S$ be the $n$-pre-kite defined by
\begin{eqnarray*} S&=& PK[n;u;\vv] ~=~ PK[n;u;v_1,\cdots,v_n]. 
  \end{eqnarray*}
For  $0 \le j \le n$, let  $S_j$ be the $j$-th facet of $S$, and let $\CCC_j$ and $\DDD_j$ be the  Cayley-Menger and the inner Cayley-Menger determinants of $S_j$. Let
\begin{eqnarray} \label{alphabeta} 
\alpha = u + v_1 + \cdots + v_n,~~\beta = u^2+v_1^2+\cdots+v_n^2.
\end{eqnarray}
Then for $1 \le j \le n$, we have
\begin{eqnarray}
\CCC_0  &=&  (-1)^n n u^{n-1}. \label{C0}\\
\CCC_j &=& (-u)^{n-3} \left[ - \alpha^2 + (n-1) \beta  - n v_j^2   + 2 \alpha v_j \right]. \label{Cj} \\
\DDD_0  &=& (-1)^{n+1} u^n (n-1).\label{D0} \\
\DDD_j  &=& (-u)^{n-2} [ (n-2)\beta -\alpha^2 + 2 \alpha u  - (n-1) u^2 - (n-1) v_j^2   \nonumber \\&&  + 2 \alpha v_j  - 2u v_j ].\label{Dj}
\end{eqnarray}
\end{thm}

\vspace{.2cm} \noindent {\it Proof.} Observing that 
\begin{eqnarray}
S_0&=& PK[n-1;u;u,\cdots,u]\\
S_j&=& PK[n-1;u;v_1,\cdots,v_{j-1},v_{j+1},\cdots,v_n] \mbox{~~if $1\le j \le n$},
\end{eqnarray}
and using Theorem \ref{CCCDDD}, we obtain
\begin{eqnarray*}
\CCC_0  &=& (-u)^{n-3}
 \left[(n-1) nu^2  -  n^2u^2 \right] \nonumber\\
 &=& (-1)^n n u^{n-1},\\
\CCC_j  
&=& (-u)^{n-3} \left[(n-1) ( \beta - v_j^2) - (\alpha - v_j)^2 \right] \\
&=& (-u)^{n-3} \left[(n-1) \beta - (n-1) v_j^2  
- \alpha^2 - v_j^2  + 2 \alpha v_j\right] \\
&=& (-u)^{n-3} \left[- \alpha^2 + (n-1) \beta  - n v_j^2    + 2 \alpha v_j\right],  \\
\DDD_0 
&=& (-u)^{n-2} \left[(n-2)u^2 (n-1) - (n-1)^2 u^2    \right] \\
&=& (-1)^{n+1} u^{n} (n-1), \\
\DDD_j 
&=& (-u)^{n-2} \left[(n-2)(\beta -v_j^2 - u^2)-(\alpha - v_j-u)^2\right] \\
&=& (-u)^{n-2} \left[ (n-2)\beta - (n-2) v_j^2 - (n-2) u^2 -\alpha^2 - v_j^2 - u^2 + 2 \alpha v_j + 2 \alpha u - 2u v_j
\right] \\
&=& (-u)^{n-2} \left[
(n-2)\beta -\alpha^2 + 2 \alpha u  - (n-1) u^2 
- (n-1) v_j^2      + 2 \alpha v_j  - 2u v_j
\right].
 \end{eqnarray*}
This completes the proof. \hfill $\Box$

\bigskip
The next theorem uses Theorem \ref{CCCDDD} to provide another derivation of the relation (\ref{REL}) mentioned earlier. 

\begin{thm} \label{Bentin} Let   $S = [A_1,\cdots,A_{n+1}]$ be a regular $n$-simplex of edge length $t_0$, and let $P$ be a point in its affine hull. Let $t_j$, $1 \le j \le n+1$, denote the distance from $P$ to the vertex $A_j$. Then
\begin{eqnarray} \label{RELREL}
(n+1) \sum_{j=0}^{n+1} t_j^4 &=& \left( \sum_{j=0}^{n+1} t_j^2 \right)^2. 
\end{eqnarray}
\end{thm}

\vspace{.1cm} \noindent {\it Proof.} 
Since the point $P$ lies in the affine hull of the regular $n$-simplex
$S = [A_1,\cdots,A_{n+1}]$, then the $(n+1)$-pre-kite (having $P$ as an apex) is degenerate, and hence has volume 0. Equivalently, its Cayley-Menger determinant is 0. Since $t_0$ is not 0, the relation (\ref{RELREL}) follows immediately from Theorem \ref{CCCDDD} and (\ref{V}).\hfill $\Box$

\bigskip We end this section by giving a new proof of Pompeiu's theorem. 
The proof uses the last part of the next lemma; the other parts of the lemma will be used  later.

\begin{lem} \label{uvh}
Let $S = [A_0, \cdots, A_n]$ be a regular $n$-simplex with center $I$, and let $u$ and $R$ be its edge length and circumradius, respectively. Let $G$ be the center of the (regular) $(n-1)$-simplex $S_0 = [A_1, \cdots, A_n]$. Then
\begin{eqnarray}
\frac{R^2}{u^2} &=& \frac{n}{2(n+1)}, \label{(i)}\\
\|A_0 - G\|  &=& \frac{(n+1)R}{n} = \sqrt{\frac{n+1}{2n}}u. \label{(ii)}
 \end{eqnarray}
If $P$  is an arbitrary point in the affine hull of $S$ with  $\|P-I\| = \rho$, then
\begin{eqnarray}
\|P-A_0\|^2 + \cdots + \|P-A_n\|^2 &=& (n+1) (\rho^2 + R^2), \label{(iii)}
\end{eqnarray}
and therefore  $P$ lies on the circumsphere of $S$ if and only if
\begin{eqnarray}
\|P-A_0\|^2 + \cdots + \|P-A_n\|^2 &=& 2(n+1)  R^2 ~=~  n u^2. \label{(iv)}
\end{eqnarray}
\end{lem}

\vspace{.3cm} \noindent {\it Proof.} For (\ref{(i)}) and (\ref{(ii)}), see  Proposition 4.6 (p. 281) of \cite{RM-ortho}. For (\ref{(iii)}),
assume, without loss of generality,  that $I$ is the origin $\OOO$. Then
	$\|A_j\| = R$ for $0 \le j \le n$, and $\|P \| = \rho$. Also $A_0 + \cdots + A_n = \OOO$, and hence $P \cdot A_0 + \cdots + P\cdot A_n = 0$. Therefore
	\begin{eqnarray*}
\|P-A_0\|^2 + \cdots + \|P-A_n\|^2 &=& (n+1)\rho^2 + (n+1) R^2,
		\end{eqnarray*}
	as desired. For (\ref{(iv)}), we use (\ref{(iii)})  and (\ref{(i)}).
	\hfill $\Box$

\begin{thm} \label{Pompeiu} Let $T$ be an equilateral triangle with circumcircle $\Gamma$, and let $P$ be an arbitrary point in its plane. Then  the distances from $P$ to the vertices of $T$ can serve as the side lengths of a triangle $T_P$. Also, $T_P$ is degenerate if and only if $P$ lies on  $\Gamma$. 
 \end{thm}
 
 \vspace{.15cm} \noindent {\it Proof.}
 Let $a$ be the side lengths of $T$, and let  $x$, $y$, and $z$ be the distances from $P$ to the vertices of $T$. By the case $n=2$ of (\ref{REL}), we have
 \begin{eqnarray} \label{g}
g : = 3(a^4+x^4+y^4+z^4) - \left(a^2+x^2+y^2+z^2\right)^2 &=& 0.
\end{eqnarray}
This polynomial $g$  simplifies into 
\begin{eqnarray*}
g &=& 2[a^4  -a^2 (x^2+y^2+z^2)]+3(x^4+y^4+z^4)- (x^2+y^2+z^2)^2\nonumber\\
&=& 2\left( a^2  - \frac{x^2+y^2+z^2}{2}\right)^2
- \frac{3}{2} (x^2+y^2+z^2)^2
+3(x^4+y^4+z^4) \nonumber\\
&=& 2\left( a^2  - \frac{x^2+y^2+z^2}{2}\right)^2
- \frac{3}{2} 
\left[2(x^2y^2+y^2z^2+z^2x^2) -  (x^4+y^4+z^4)\right]. \nonumber\\
\end{eqnarray*}
Since $g=0$, it follows that 
\begin{eqnarray}
h:= 2(x^2y^2+y^2z^2+z^2x^2) -  (x^4+y^4+z^4) &\ge& 0, \label{ti}
\end{eqnarray}
with equality if and only if 
\begin{eqnarray} \label{equality}
a^2  &=&  \frac{x^2+y^2+z^2}{2}.
\end{eqnarray}
By the last part of  Lemma \ref{uvh}, this is equivalent to saying that $P$ is on $\Gamma$. Thus we assume that $h > 0$, i.e.,
\begin{eqnarray} \label{ti-2}
(x+y+z)(-x+y+z)(x-y+z)(x+y-z)  &>& 0.
\end{eqnarray}
Since $x+y+z > 0$,  $h >0$ is equivalent to
\begin{eqnarray} \label{ti-3}
(-x+y+z)(x-y+z)(x+y-z)  &>& 0.
\end{eqnarray}
Since the sum of any two of the terms $-x+y+z$, $x-y+z$, and $x+y-z$  is non-negative, it follows that at most one of these terms
is negative. Thus if one of them is negative, then the other two are non-negative, contradicting (\ref{ti-3}).
Therefore the three terms are non-negative, i.e., positive. Therefore 
\begin{eqnarray} \label{333}
y+z >  x,~~z+x > y, ~~ x+y > z,
\end{eqnarray}
proving that $x$, $y$, and $z$ can serve as the side lengths of a (non-degenerate) triangle.

Thus if $h=0$, $P$ lies on $\Gamma$, and $x$, $y$, and $z$ form the side lengths of a degenerate triangle; if $h > 0$, $P$ does not lie on $\Gamma$, and $x$, $y$, and $z$ form the side lengths of a non-degenerate triangle. This is what we were to prove.
\hfill $\Box$


\section{Two-apexed pre-kites and the limitations on their edge lengths} \label{Jaws} 
We define a {\it two-apexed $n$-pre-kite} to be an $n$-simplex $S = [A_0,\cdots,A_n]$, $n \ge 3$,  in which two of its facets are regular  $(n-1)$-simplices.   Notice that a two-apexed  $n$-pre-kite is nothing but an $n$-pre-kite with two apexes. It will also be proved in Lemma \ref{1-2}
that a non-regular  $n$-pre-kite cannot have more than 2 apexes. 

The main theorem in this section, namely Corollary \ref{cor},  gives necessary and sufficient conditions on given  positive numbers  so that they can serve as the edge lengths of a two-apexed $n$-pre-kite.

\bigskip
\begin{lem} \label{1-2}
A non-regular  $n$-pre-kite, $n \ge 3$,  can have at most two apexes.
\end{lem}

\vspace{.15cm} \noindent {\it Proof.} Let $S = [A_0, \cdots, A_n]$ be an $n$-pre-kite, and suppose that $A_0$, $A_1$,  and $A_2$ are three apexes. We are to prove that $S$ is regular. Since  $A_2A_3$
is an edge in both facets $S_0$ and $S_1$, and since these facets are regular, it follows that the facets $S_0$ and $S_1$ have the same edge length. Similary, we show that 
the facets $S_0$ and $S_2$ have the same edge length.
But every edge in $S$ is an edge in one of the facets  $S_0$, $S_1$, and $S_2$. Thus all edges of $S$ have the same length. \hfill $\Box$

\begin{thm} \label{limitations}
Let $S = [A_0, \cdots, A_n]$ be a regular $n$-simplex  that lies in $\RR^m$  for some $m \ge n+1$, and let its side length be $u$ and its circumradius be $R$. Let $Q_0$ be the reflection of $A_0$ about the affine hull $H_0$ of the facet $S_0 = [A_1,\cdots,A_n]$, and let $\Omega$ be the set of all points $Q$ in $\RR^m$ for which the  $n$-simplex $[Q, A_1, \cdots, A_n]$ is regular. 
Then 
\begin{eqnarray}
\{\|Q - A_0\| : Q \in \Omega\} &=& \left[ 0, \sqrt{\frac{2(n+1)}{n}} u\right] \label{lim}\\
&=& \left[ 0, \frac{2(n+1)}{n}  R  \right], \label{limm}
\end{eqnarray}
 with the extreme values taken at $Q=A_0$ and at $Q=Q_0$.
\end{thm}

\vspace{.15cm} \noindent {\it Proof.} Without loss of generality, we assume that the center of $S$ lies at the origin $\OOO$ of $\RR^m$.

Let $\CCC$ be the collection of all $(n+1)$-dimensional subspaces of $\RR^m$ that contain $S$, and for any $V \in \CCC$, let $L_V = \{\|Q - A_0\| : Q \in \Omega \cap V\}$. Let $LHS$ and $RHS$ stand for the left and right hand sides of (\ref{lim}), respectively. If we could prove that $L_V = RHS$, then we will be done. This is because every point in $\RR^m$ belongs to some $V \in \CCC$,  and hence $LHS$ is the union of all $L_V$, where $V$ ranges in $\CCC$, and since $RHS$ depends on $S$ only. Thus we take any $V \in \CCC$, and we are to prove that $L_V = RHS$. In other words, $$\mbox{we assume that $m = n+1$, and we set $\RR^m = V$}.$$


Let $G$ be the center of $S_0$, and let $H_0$ be the affine hull of $S_0$. Let $W = \{ P \in V : (P-G) \perp H_0\}$. Thus $W$ is the shifted orthogonal complement of $H_0$, namely $W - G = (H_0 - G)^{\perp}$. Thus $\dim W = \dim V - \dim H_0 = (n+1)-(n-1) = 2,$ i.e., $W$ is a plane. We have used the facts that 
if $A$ is a closed subspace of a Hilbert space $B$, then $B = A \oplus A^{\perp}$ (\cite[Theorem 3.3-4, p.~146]{Kry}), and that every finite dimensional subspace  of a normed space is closed (\cite[Theorem 2.4-5, p.~74]{Kry}).
Since $G$ is the center of the facet $S_0$ of the regular $n$-simplex $S$, it follows that $A_0 - G$ is an altitude of $S$, i.e., $(A_0-G) \cdot H_0 = 0$, and hence  $A_0 \in W$. Similarly, $Q_0 \in W$.

Let $\Gamma$ be the circle in $W$ centered at $G$ and passing through $A_0$ (and $Q_0$).  We claim that $\Omega = \Gamma$. To see this, let $P \in \Gamma$, and let  $1 \le i \le n$. Then $P \in W$ and hence $(G-P) \perp H_0$. Therefore
\begin{eqnarray*}
\|P-A_i\|^2 &=& \|P- G\|^2 + \| G-A_i\|^2, \mbox{~by Pythagoras' theorem}\\
&=& \|A_0 - G\|^2 + \| G- A_i\|^2, \mbox{~because $P \in \Gamma$}\\
&=& \|A_0- A_i\|^2, \mbox{~by Pythagoras' theorem}.
\end{eqnarray*}
This shows that $[P,A_1,\cdots,A_n]$ is regular, and therefore $P \in \Omega$. Conversely, let $P \in \Omega$. Thus $T=[P, A_1,\cdots,A_n]$ is regular. Since $G$ is the center of the facet $S_0$ of $T$, and since $T$ is regular,
it follows that
$(P-G)\perp H_0$, and hence $P \in W$. In particular, $A_0, Q_0 \in W$. 
Therefore
\begin{eqnarray*}
\|P-G\|^2 &=& \|P- A_1\|^2 - \| A_1 -G\|^2, \mbox{~by Pythagoras' theorem}\\
&=& \|A_2- A_1\|^2 - \| A_1 -G\|^2, \mbox{~because $[P,A_1,\cdots,A_n]$ is regular}\\
&=& \|A_0- A_1\|^2 - \| A_1 -G\|^2, \mbox{~because $[A_0,A_1,\cdots,A_n]$ is regular}\\
&=& \|A_0- G\|^2, \mbox{~by Pythagoras' theorem}.
\end{eqnarray*}
Therefore $P \in \Gamma$. Thus we have shown that $\Gamma = \Omega$.

Since $S=[A_0,\cdots,A_n]$ is regular and since $G$ is the center of its facet  
$S_0 = [A_1,\cdots,A_n]$, it follows that $A_0G$ is perpendicular to $H_0$. Also $A_0Q_0$ is perpendicular to $H_0$. Therefore  $A_0$, $G$, and $Q_0$ are collinear. Since $A_0$  and $Q_0$ lie on $\Gamma$, and since $G$ is the center of $\Gamma$, it follows that $A_0Q_0$ is a diameter of $\Gamma$ (with midpoint $G$). Therefore, as $Q$ moves on $\Gamma$, $\|Q-A_0\|$ takes all values between $0$ and $\|Q_0-A_0\|$. Thus our proof will be complete if we prove that  $$\|Q_0-A_0\| = \frac{2(n+1)}{n}   R ~=~ \sqrt{\frac{2(n+1)}{n}} u.$$
But this follows immediately from the formula for $\| A_0 - G\|$ given in  (\ref{(ii)}) and the fact that
\begin{eqnarray*}
\|A_0 - Q_0\| &=& 2 \|A_0 - G\|.
\end{eqnarray*}
This completes the proof. 
\hfill $\Box$

\bigskip The following corollary will be used in a later section.

\begin{cor} \label{cor}  There exists an $n$-simplex having one edge of length $v > 0$ and having all the remaining edges of lengths $u > 0$ if and only if 
\begin{eqnarray}\label{uv}
0 < \frac{v}{u} < \sqrt{\frac{2n}{n-1}}.
\end{eqnarray}
In other words, there exists a two-apexed  $n$-pre-kite $PK[n;u;v_1,\cdots, v_n]$, with $v_1=\cdots,v_{n-1}=u$ and $v_n=v$ if and only if $u$ and $v$ satisfy 
(\ref{uv}).
\end{cor}

\vspace{.2cm} \noindent {\it Proof.} The previous theorem shows that
there exists an $(n+1)$-simplex having one edge of length $v > 0$ and having all the remaining edges of lengths $u > 0$ if and only if 
\begin{eqnarray}
0 < \frac{v}{u} < \sqrt{\frac{2(n+1)}{n}}.
\end{eqnarray}
The desired result is obtained by replacing $n+1$ by $n$. \hfill $\Box$

\section{Coincidence of two of the classical centers of a pre-kite} \label{coincidence}
The classical centers of an $n$-simplex $S=[A_0,\cdots,A_n]$ refer to the circumcenter, the incenter, and the centroid
 of $S$. The {\it circumcenter}  $\QQQ = \QQQ (S)$ of $S$ is the center of the $(n-1)$-sphere that passes through the vertices of $S$. The {\it incenter}  $\III = \III (S)$ of $S$ is the center of the $(n-1)$-sphere that touches the facets internally, i.e., at points that lie in the convex hulls  of the facets. The {\it centroid}  $\GGG = \GGG (S)$ of $S$ is defined inductively to be the intersection of the medians  of $S$, where a {\it median} of $S$ is the line segment joining a vertex of $S$ to the centroid of the opposite facet. It is also defined by the simple formula
 $$\GGG (S) = \frac{A_0 + \cdots + A_n}{n+1}.$$

Theorem \ref{O=G} proves that if the circumcenter and the centroid of an $n$-pre-kite $S$, $n \ge 2$, coincide, then it is regular. 
Theorem \ref{O=I} proves that if the circumcenter and the incenter of an $n$-pre-kite $S$, $n \ge 2$, coincide, then it is regular. 
Theorem \ref{I=G} proves that if the incenter and the centroid of an $n$-pre-kite $S$ coincide, and if $n \le 5$, then it is regular, and exhibits examples of non-regular $n$-pre-kite $S$, $n \ge 6$, in which  the incenter and the centroid coincide.

The main  tools in proving Theorems  \ref{O=G}, \ref{O=I}, and \ref{I=G}
are the formulas established in Theorem \ref{CjDj}, together with 
the following theorem, proved in \cite[Theorem 3.2, p.~496]{BAG-1}. We recall that an $n$-simplex is said to be {\it equiareal} if its facets have equal  volumes, i.e., equal  $(n-1)$-dimensional Lebesgue measures. 
It is said to have {\it well-distributed} edge lengths if its facets have equal variance, i.e., if the sum of squares of the edge lengths of a facet is the same for all facets. It is said to be  {\it equiradial}  if its facets have equal  circumradii.

\begin{thm} \label{BAG-1}
Let $S = [A_0,\cdots,A_n]$ be an $n$-simplex. 
Then
\begin{enumerate}
	\item[(i)] The circumcenter $\QQQ$ and the centroid $\GGG$ of $S$ coincide if and only if $S$ has {\it well-distributed} edge lengths.
	\item[(ii)] The circumcenter  $\QQQ$  and the incenter $\III$ of $S$ coincide if and only if $\QQQ$ is interior and $S$ is equiradial.
	\item[(iii)] The centroid $\GGG$  and the incenter $\III$ of $S$ coincide if and only if $S$ is equiareal.
	\item[(iv)] The centroid $\GGG$, the circumcenter $\QQQ$, and the incenter $\III$ of $S$ coincide if and only if two of the conditions
	\begin{itemize}
	\item[(a)] $S$ has well-distributed edge lengths,
	\item[(b)] $S$ is equiradial, 
		\item[(c)] $S$ is equiareal
		\end{itemize}
	hold. When this happens, the third condition also holds.
	\end{enumerate}
\end{thm}

\begin{thm} \label{O=G}
Let $S = [A_0,\cdots,A_n]$ be an $n$-pre-kite, and suppose that $n \ge 2$. If  the circumcenter  $\QQQ$ and the centroid  $\GGG$ of $S$ coincide,
then $S$ is regular. \end{thm}

\vspace{.1cm}\noindent {\it Proof.} 
Suppose that  the circumcenter  $\QQQ$ and the centroid  $\GGG$ of $S$ coincide. By Theorem \ref{BAG-1}, $S$ has  well-distributed edge lengths. Then the sum $L_j$, $0 \le j \le n$, of squares of the lengths of the edges of the facet $S_j$ does not depend on $j$. Let $L$ be
the sum of squares of the lengths of the edges of $S$, and let 
$M_j$, $0 \le j \le n$, be the sum of squares of the lengths of the edges that emanate from $A_j$. Since $M_j = L - L_j$, it follows that $M_j$ does not depend on $j$. 
Assuming that $A_0$ is the apex of the pre-kite $S$, it follows from (\ref{v}) that
\begin{eqnarray}
M_0 &=& v_1^2 + \cdots + v_n^2, \label{00}\\
M_j &=& v_j^2 + (n-1) u^2 \mbox{~~for $0 \le j \le n$.} \label{jj}
\end{eqnarray}
It follows from (\ref{jj}) that
$v_1 = \cdots = v_n$. If $v$ is the common value of 
$v_1, \cdots,  v_n$, then  it follows by subtracting  (\ref{00}) from (\ref{jj}) and using that their left hand sides are equal that $(n-1)u^2 = (n-1) v^2$, and hence $v=u$. This shows that $S$ is regular, and ends the proof. \hfill $\Box$

\begin{thm} \label{O=I}
Let $S = [A_0,\cdots,A_n]$ be an $n$-pre-kite, and suppose that $n \ge 2$. If  the circumcenter  $\QQQ$ and the incenter $\III$ of $S$ coincide,
then $S$ is regular. \end{thm}

\vspace{.1cm}\noindent {\it Proof.} 
Suppose that  the circumcenter  $\QQQ$ and the incenter  $\III$ of $S$ coincide. By Theorem \ref{BAG-1}, $S$ is equiradial. Let $R_j$, $0 \le j \le n$, be the circumradius of the $j$-th facet. Let $1 \le j \le n$. Using (\ref{R}), we see that $R_0 = R_j$ if and only if 
\begin{eqnarray} \label{CD01}
\CCC_0 \DDD_j &=& \CCC_j \DDD_0. 
\end{eqnarray}
By Theorem \ref{CjDj}, (\ref{CD01})  simplifies into
\begin{eqnarray} \label{j-ind}
2(\alpha - nu)v_j&=&\alpha^2 + \beta - 2n\alpha u+n(n-1)u^2.
\end{eqnarray}
If $\alpha = nu$, then it follows from (\ref{j-ind}) that $\beta = nu^2$. By (\ref{CCC-PK-2}), the Cayley-Menger determinant $\CCC (PK[n;u;\vv])$ of $S$ is 0. Thus $S$ is degenerate, which we discard. Therefore $\alpha \neq nu$. It now follows from (\ref{j-ind})  that $v_j$ does not depend on $j$. Therefore  $v_1=\cdots=v_n$. This means that $S$ is a $n$-kite. By Lemma 4.5 of \cite{RM-ortho},  
$S$ is regular. 
 \hfill $\Box$

\begin{thm}\label{I=G}
Let $S = [A_0,\cdots,A_n]$, $n \ge 2$,  be an $n$-pre-kite.
Suppose that  the incenter $\III$  and the centroid $\GGG$ of $S$ coincide, i.e., $S$ is equiareal.

If $n\le 5$, then $S$ is regular.

If $n \ge 6$, then $S$ is not necessarily regular; i.e. there exist non-regular  $n$-pre-kites, in fact two-apexed $n$-pre-kites,  in which the incenter and centroid coincide. 
 \end{thm}

\vspace{.15cm}\noindent {\it Proof.} Suppose that $S=[A_0,\cdots,A_n]$ is a non-regular equiareal  $n$-pre-kite with apex $A_0$, say
$$S = PK [n; u; v_1,\cdots,v_n].$$

Let  $\VVV_j$ be the volume of the $j$-th facet $S_j$ of $S$, 
and let $\CCC_j$ be the Cayley-Menger determinant of $S_j$, as defined in (\ref{CCC}). By Theorem \ref{CjDj}, 
$$
\CCC_j = \left\{ \begin{array}{ll} 
(-u)^{n-3} \left[ - \alpha^2 + (n-1) \beta  - n v_j^2   + 2 \alpha v_j \right] 
& \mbox{if $1\le j\le n$},\\
\vspace{.15cm}
(-1)^n n u^{n-1} & \mbox{if $j = 0$}.
\end{array} \right.
$$
By (\ref{V}),  $\VVV_i = \VVV_j \ifff \CCC_i = \CCC_j$.

We prove first that  $v_1, \cdots, v_n$ cannot be all equal.   In fact, if $v_1=\cdots=v_n (= x, \mbox{~say})$, then the condition $\VVV_1 = \VVV_0$ 
 yields $u=x$ as follows:
\begin{eqnarray*}
\VVV_1 = \VVV_0 &\ifff& nu^2 - \alpha^2 + (n-1)\beta  - nx^2 + 2 \alpha x  = 0\\
&\ifff& 2(n-1)u^2     + 2ux (1-n)  = 0\\
&\ifff& 2(n-1)u (u-x) = 0\\
&\ifff& u = x.
\end{eqnarray*}
Then $S$ is regular, contradicting the assumptions.

Thus we assume that  $v_1,\cdots,v_n$ are not all equal.

Next, we prove  that there do not exist three distinct indices $i, j, k$ in $\{1, 2, \cdots, n\}$ such that $v_i$, $v_j$, and $v_k$ are pairwise different. This is because the existence of such indices  contradicts the assumption $\VVV_i = \VVV_j = \VVV_k$. In fact, \begin{eqnarray*}
&&\VVV_i = \VVV_j = \VVV_k ~\ifff~  \CCC_i = \CCC_j = \CCC_k\\
&\ifff& - n v_i^2   + 2 \alpha v_i = - n v_j^2   + 2 \alpha v_j = - n v_k^2   + 2 \alpha v_k \\
&\ifff& (v_i - v_j)  \left[- n (v_i + v_j)    + 2 \alpha\right] =
(v_j - v_k) \left[- n (v_j + v_k)    + 2 \alpha  \right] = 0
\\
&\ifff& - n (v_i + v_j) + 2 \alpha=- n (v_j + v_k)    + 2 \alpha = 0\\ 
&\riff& - n (v_i + v_j)  = -n (v_j + v_k)
\\ &\riff& v_i = v_k,
\end{eqnarray*}
a contradiction. Therefore no three of the numbers $v_1,\cdots,v_n$ are pairwise distinct. 

Thus there  are two different numbers $x$ and $y$ and an index $t \in \{1, \cdots, n\}$ such that
\begin{eqnarray}
\mbox{$v_j = x$ if $1 \le j \le t$, and $v_j = y$ if $t <  j \le n$}. \end{eqnarray}
Let $s = n-t$. We may clearly assume that $t \ge s$. Thus $S$ is of the form 
\begin{eqnarray}
S &=& PK[n;u;\overbrace{x,\cdots,x}^t,\overbrace{y,\cdots,y}^s], \end{eqnarray}
where $x$ is repeated $t$ times and $y$ is repeated $s = n-t$ times.
Let $i$ and $j$ be such that $v_i = x$ and $v_j = y$.
Then 
\begin{eqnarray*}
\VVV_i =  \VVV_j &\ifff& - n x^2   + 2 \alpha x = - n y^2   + 2 \alpha y\\
&\ifff& (x-y)  \left[- n (x+y)    + 2 \alpha  \right]=0\\
&\ifff& - n (x+y)    + 2 \alpha  =0\\
&\ifff& - n (x+y)    + 2 (tx + sy)   + 2 u =0\\
&\ifff& (-n + 2t) x+ (-n+2s) y +2u   =0\\
&\ifff& (t-s) x+ (s-t) y +2u =0\\
&\ifff& (t-s)(y-x) = 2u.
\end{eqnarray*}
Also, 
\begin{eqnarray*}
\VVV_i =   \VVV_0  &\ifff& - \alpha^2 + (n-1)\beta - nx^2 + 2 \alpha x  = -nu^2\\
&\ifff& nu^2 - \alpha^2 + (n-1)\beta  - nx^2 + 2 \alpha x = 0.
\end{eqnarray*}
Therefore
\begin{eqnarray}
\mbox{$S$ is equiareal}  &\ifff& \VVV_i =  \VVV_0 \mbox{~and~} \VVV_i= \VVV_j\nonumber \\
&\ifff& (i)~ nu^2 - \alpha^2 + (n-1)\beta  - nx^2 + 2 \alpha x  = 0, \mbox{~and~}\nonumber \\
&&(ii) ~(t-s)(y-x) = 2u. \label{i-ii}
\end{eqnarray}
Thus equiareality of $S$ is equivalent to fulfilment of (i) and (ii) of (\ref{i-ii}).  
Since $u \ne 0$, these imply that 
\begin{eqnarray}
t &\ne& s. \label{ts} 
\end{eqnarray}

Let us first treat the case $$s=1,~~t=n-1.$$ In this case,
\begin{eqnarray*}
&&\mbox{$S$ is equiareal}\\
&\ifff& (i) ~2nu^2 - \alpha^2 + (n-1)\beta  - nx^2 + 2 \alpha x  = 0, \mbox{~and~}\\
&&(ii)~ (n-2)(y-x) = 2u.
\end{eqnarray*}
Plugging 
$2u = (n-2)(y-x)$ in (i), we obtain
$(n-1) (n-2) (x-y) (x n - y n + 2 y)=0$, i.e., $x n - y n + 2 y=0$. Solving this with $2u = (n-2)(y-x)$, we obtain $x=u$. Thus one of the edge lengths of $S$ is
$$\frac{un}{n-2},$$ and each other  edge length is $u$.
In view of Corollary \ref{cor}, 
\begin{eqnarray*}
\mbox{such an $n$-simplex exists} &\ifff& \frac{n}{n-2} < \sqrt{\frac{2n}{n-1}}\\
&\ifff&      n^2 (n-1)   <  2n (n-2)^2 \\
&\ifff&          n^2 - 7n + 8  > 0\\
    &\ifff&          n > \frac{7 + \sqrt{17}}{2} \approx 5.6\\
    &\ifff&          n \ge 6.
    \end{eqnarray*}
Thus if $n \ge 6$, there are non-regular equiareal $n$-pre-kites. These can even be chosen to be of the form $PK[n;u;v_1,\cdots,v_n]$ with $v_1=\cdots=v_{n-1}=u$, i.e., a two-apexed $n$-pre-kite.
If $n \le 5$, then a non-regular equiareal $n$-pre-kite must have 
\begin{eqnarray}
t, s  &\ge& 2. \label{ge2} 
\end{eqnarray}
In view of (\ref{ts}), this is  possible only if $(n,t,s) = (5,3,2)$. We show now that this cannot happen  either. In fact, this assumption would imply that
$$\alpha = u + 3x + 2y,~~\beta = u^2 + 3x^2 + 2y^2,~~2u = y- x.$$
Substituing these values of $\alpha$, $\beta$, and $u$ in (ii) and factorizing,  we obtain $$(y-x)(y-2x) = 0.$$
Since $y \ne x$, it follows that $y=2x$ and $u = x/2$. Thus our $5$-pre-kite is of the form $PK[5;u;2u,2u,2u,4u,4u]$. Using (\ref{CCC-PK-2}), we see that the Cayley-Menger determinant of such a pre-kite is
$$(-u)^{3} [5(u^2+3(4u^2)+2(16u^2)) - (u + 3(2u)+2(4u))^2] = (-u)^{3} [(5)(45u^2) -  (15u)^2)] = 0.$$
Thus this pre-kite, if it exists, is degenerate, which we reject.

This completes the proof. \hfill $\Box$

\begin{que}
{\rm The {\it Fermat-Torricelli point}  $\FFF = \FFF (S)$  of an $n$-simplex $S$ is the point whose distances from the vertices of $S$ have a minimal sum. It is often thought of as a semi-classical (or even a classical) center. Thus it is natural to investigate  the degrees of regularity implied by the coincidences
$\FFF=\GGG$, $\FFF=\III$, and $\FFF=\QQQ$. In this regard, we recall  Theorem 3.1, p. 496, of \cite{BAG-1}. This states that if any two of the three centers $\FFF$, $\QQQ$, and  $\GGG$ coincide, then all the three coincide. Thus each of the coincidences $\FFF=\GGG$ and $\FFF=\QQQ$ implies that $\QQQ=\GGG$, and hence regularity (by Theorem \ref{O=G} above). This leaves us with the question about the degree of regularity implied by the coincidence $\FFF=\III$. We leave this open.}
\end{que}

\section{Pre-kites in the four special families of simplices} In this section, we shall see how the new family  of $n$-pre-kites is related to the four known special families  of orthocentric, circumscriptible, isodynamic, and tetra-isogonic $n$-simplices. 

We recall that an $n$-simplex $S = [A_1,\cdots,A_{n+1}]$, $n \ge 2$,  is 
said to be {\it orthocentric}   if  the altitudes of $S$ are concurrent. 
It is said to be {\it circumscriptible} (or {\it edge-incentric} or {\it balloon}) if there is an $(n-1)$-sphere that touches all its  edges internally.
It is said to be {\it isodynamic} if  the incentral cevians are concurrent. Here, an {\it incentral cevian} is the cevian that joins a vertex and the incenter of the opposite facet. It is said to be  {\it tetra-isogonic} if every four vertices of $S$ form an isogonic tetrahedron, i.e., a tetrahedron whose inspherical cevians are concurrent. Here, an {\it inspherical cevian} is the cevian that joins a vertex to the point where the insphere touches the opposite face. Other characterizations appear in \cite{RM-HHM}.

Let us denote the  families  of orthocentric, circumscriptible, isodynamic, and tetra-isogonic $n$-simplices and the families of  $n$-kites  and  $n$-pre-kites  by $F_o$, $F_c$, $F_d$,  $F_g$, $F_k$, and $F_p$, respectively. 

It is proved in \cite{kites} that the intersection of any two of the families \begin{eqnarray}\label{4four}
F_o, ~F_c, ~F_d, ~\mbox{and~} F_g
\end{eqnarray}
is the family $F_k$. In this section, we prove that this still holds if we enlarge the list in (\ref{4four}) to include our new family $F_p$. We prove this in Theorem \ref{5five}. The proof is a consequence of the following theorem, which is taken from \cite{kites}.

\begin{thm} \label{HHM}
Let  $S = [A_1,\cdots,A_{d+1}]$, $d \ge 2$, be a $d$-simplex. Then
\begin{eqnarray*} 
\mbox{$S$ is orthocentric} &\ifff&\mbox{there exist $\beta_1, \cdots, \beta_{d+1} \in \RR$ such that} \nonumber \\ 	
&&\mbox{$\|A_i-A_j\|^2=\beta_i + \beta_j$
for $1 \le i <  j \le d+1$}\\
\mbox{$S$ is circumscriptible} &\ifff& \mbox{there exist $\beta_1, \cdots, \beta_{d+1} > 0$ 		such that} \nonumber \\ 	&&\mbox{$\|A_i-A_j\|=\beta_i + \beta_j$ for $1 \le i <  j \le d+1$}\\
\mbox{$S$ is isodynamic} 	&\ifff&\mbox{there exist $\beta_1, \cdots, \beta_{d+1} > 0$ 		such that} \nonumber \\ 	&&\mbox{$\|A_i-A_j\|^2 =\beta_i  \beta_j$ for $1 \le i <  j \le d+1$}\\
\mbox{$S$ is tetra-isogonic} 	&\ifff&\mbox{there exist $\beta_1, \cdots, \beta_{n+1} > 0$ 		such that} \nonumber \\ 	&&\mbox{$\|A_i-A_j\|^2=\beta_i^2 + \beta_i \beta_j +   \beta_j^2$ for $1 \le i <  j \le d+1$}.
\end{eqnarray*}
Moreover, the numbers $\beta_i$, $1 \le i \le d+1$, appearing in the four equations are unique.
\end{thm}

\begin{thm} \label{5five}
For $n \ge 3$, the intersection of any two of the five families
\begin{eqnarray}\label{4}
F_o, ~F_c, ~F_d, ~F_g, ~\mbox{and~} F_p
\end{eqnarray}
is the family $F_k$.
\end{thm}

\vspace{.2cm} \noindent{\it Proof.}
In view of the fact, proved in \cite{kites}, that
the intersection of any two of the four families $F_o$, $F_c$, $F_d$,  and $F_g$ is the family $F_k$, it remains to show that
\begin{eqnarray}\label{55five}
F_o \cap F_p = F_c \cap F_p =  F_d \cap F_p =  F_g \cap F_p = F_k.
\end{eqnarray}
For this, Theorem \ref{HHM} is very useful. Since the proofs of these statements are similar, we find it sufficient to prove the last statement only, i.e.,
\begin{eqnarray}\label{FgFp}
F_g \cap F_p = F_k.
\end{eqnarray}
Thus let $S = [A_0,\cdots,A_n]$, $n \ge 3$,  be the $n$-pre-kite
$PK [n;u;v_1,\cdots,v_n]$ with apex $A_0$, and suppose that $S$ is 
tetra-isogonic. By Theorem \ref{HHM}, there
exist $\beta_0, \cdots,  \beta_{n} > 0$ 		such that 
\begin{eqnarray}
\|A_i-A_j\|^2=\beta_i^2 + \beta_i \beta_j +   \beta_j^2 \mbox{~~ for $0 \le i <  j \le n$}.
\end{eqnarray}
By the definition of $PK[n;u;v_1,\cdots,v_n]$, we see that
\begin{eqnarray}
v_j &=& \beta_0^2 + \beta_0 \beta_j +   \beta_j^2 \mbox{~~ for $1 \le   j \le n$},\\
u &=& \beta_i^2 + \beta_i \beta_j +   \beta_j^2 \mbox{~~ for $1 \le i <  j \le n$}.
\end{eqnarray}
It follows from the second equation that
\begin{eqnarray}
u= \beta_1^2 + \beta_1 \beta_i +   \beta_i^2 
&=&
\beta_1^2 + \beta_1 \beta_j +   \beta_j^2 
\mbox{~~ for $2 \le i <   j \le n$}.
\end{eqnarray}
Therefore $(\beta_i - \beta_j) (\beta_1 + \beta_i + \beta_j) = 0$
for $2 \le i <   j \le n$.
Since $\beta_1 + \beta_i + \beta_j > 0$, it follows that
$\beta_i = \beta_j$ for $2 \le i <   j \le n$. By symmetry, we conclude that
$\beta_i = \beta_j$ for $1 \le i <   j \le n$.  Letting $\beta$ be the common value of $\beta_1, \cdots, \beta_n$, we see that
\begin{eqnarray}
v_j  &=&  \beta_0^2 + \beta_0 \beta +   \beta^2 \mbox{~~ for $1 \le   j \le n$}.
\end{eqnarray}
This shows that $S$ is an $n$-kite. \hfill $\Box$

\end{document}